\documentstyle[12pt]{article}
\begin{document}

{\bf Space-Time Foam Dense Singularities \\ and de Rham Cohomology} \\
\\
Anastasios Mallios \\
Institute of Mathematics, University of Athens \\
Panepistimiopolis, 15784 Athens, Greece \\
AMALLIOS@ATLAS.UOA.GR \\
\\
Elemer E Rosinger \\
Department of Mathematics and Applied Mathematics \\
University of Pretoria, 0002 Pretoria, South Africa \\
EEROSINGER@HOTMAIL.COM \\

\begin{abstract}

In an earlier paper of the authors it was shown that the sheaf theoretically based recently developed abstract
differential geometry of the first author can in an easy and natural manner incorporate singularities on arbitrary closed
nowhere dense  sets in Euclidean spaces, singularities which therefore can have arbitrary large positive Lebesgue
measure. As also shown, one can construct in such a singular context a de Rham cohomology, as well as a short
exponential sequence, both of which are fundamental in differential geometry. In this paper, these results are
significantly strengthened, motivated by the so called {\it space-time foam structures in general relativity}, where
singularities can be {\it dense}. In fact, this time one can deal with singularities on arbitrary sets, provided that their
complementaries are dense, as well. In particular, the cardinal of the set of singularities can be {\it larger} than that of
the nonsingular points.

\end{abstract}

\bigskip
Note : This paper is an augmented version of the paper with the same title, published in Acta Applicandae
Mathematicae 67(1):59-89,2001, and it is posted here with the kind permission of Kluwer Academic Publishers. \\

\bigskip
'We do not possess any method at all to derive \\
systematically solutions that are free of singularities ... '

\medskip
Albert Einstein \\
{\it The Meaning of Relativity} \\
Princeton Univ. Press, 1956, p. 165 \\

\bigskip
'Sensible mathematics involves neglecting a quantity \\
when it turns out to be small - not neglecting it just \\
because it is infinitely great and you do not want it.'

\medskip
P.A.M. Dirac \\
{\it Directions in Physics} \\
H. Hora, J.R. Shepanski, Eds., J. Wiley, 1978, p. 36 \\ \\

\bigskip
{\bf 1.~ Basics of Abstract Differential Geometry} \\

\bigskip
{\bf 1.1. Introduction}

\bigskip
There is a longer tradition in enlarging the framework of classical smooth differential geometry
in such a way that singularities and various nonsmooth entities need no longer be treated as troublesome exceptions
and breakdowns in the otherwise smooth mathematical machinery, but instead, are included in it from the very
beginning, and thus can be dealt with in the same unified way, see Souriau [1,2], Geroch [1,2], Geroch \& Traschen,
Kirillov [1,2], Mostov, Blattner, Heller [1-3], Heller \& Sasin [1-3], Gruszczak \& Heller, Heller \& Mularzynski, Sasin
[1,2], Sikorski, Brylinski, Hawking \& Penrose, Penrose et.al., Mallios [1-8], Vassiliou [1-5]. \\
\\
The most far reaching approach in this respect is that recently published in Mallios [1], where instead of smooth
functions as structure coefficients, one starts with a sheaf of algebras on an arbitrary underlying topological space.
Further, one deals with a sequence of sheaves of modules, interrelated with suitable so called differentials, that is,
sheaf morphisms which satisfy a Leibniz type rule of product derivative. In this way, one obtains the differential
complex and one can recover much of the essence of classical smooth differential geometry, including de Rham
cohomology, short exact exponential sequences, characteristic classes ( \`{a} la Chern-Weil ), etc. \\
\\
In a way, the abstract differential geometry in Mallios [1] recalls what happened in general topolgy at the beginning of
the XX-th century, when it was found out that metric concepts, although originated the develpoment of topology, were in
fact not necessary, and instead, one could start with the abstract axioms of open, or equivalently, closed sets, as set
up by Kuratowski, and still recover much of the important aspects of topology. Similarly, in Mallios [1] it is shown that
Calculus, and hence, smooth functions, are in fact not necessary in developing differential geometry. Instead, suitable
sheaves of algebras of functions on rather arbitrary topological spaces can be used. And as shown in Mallios \&
Rosinger, one can go much further, by using sheaves of algebras of generalized functions. \\
\\
Earlier, in Heller [1[, see also Heller [2,3], Heller \& Sasin [1-3], Gruszczak \& Heller, Heller \& Mularzynski, Sasin [1,2],
Sikorski, it was shown that in some of such enlarged - but rather more particular - frameworks of differential geometry,
one can capture singularities in the underlying topological space which are concentrated on one single {\it fixed} closed
nowhere dense set, a set which however must be on the {\it boundary} of the space. \\
\\
In Mallios \& Rosinger, with the use of the abstract theory in Mallios [1], and of the nonlinear algebraic theory of
generalized functions in Rosinger [1-10], see also Rosinger [11-18], Oberguggenberger \& Rosinger, Rosinger \&
Rudolph, Rosinger \& Walus [1,2],  it was recently shown that the singularities can now be concentrated on
{\it arbitrary} closed nowhere dense subsets of the underlying topological space, and thus need no longer be in a
set fixed in the boundary. \\
By the way, we should recall that closed nowhere dense subsets in Euclidean spaces can have {\it arbitrary large}
positive Lebesgue measure, Oxtoby. \\
\\
In this paper, we go far further, by using a recent significant extension of the nonlinear algebraic theory of
generalized functions, see Rosinger [14-18]. Indeed, this time the singularities can be on {\it arbitrary}, including {\it
dense} subsets of the underlying topological space, provided that the complementary of the singularities, that is,
the set of nonsingular points, is itself still a {\it dense} subset. In the case of many topological spaces of interest,
in particular, Euclidean spaces, or finite dimensional smooth manifolds, this means that the singularities can be dense
and also have a cardinal {\it larger} than that of the nonsingular points. Indeed, the singularities can have the cardinal
of the continuum, while their complementary, that is, the nonsingular points, need only be countable and dense. \\
As an example, in the case of the real line ${\bf R}$, the singularities can be all the irrational points, while the
nonsingular points can be reduced to the rational ones only. \\
Interest in {\it dense} singularities arises, among others,  from the study of the so called {\it space-time foam} in
general relativity, see for instance Heller [4], or even Heller \& Sasin [3]. \\
\\
For the sake of brevity, here we shall only set up the differential complex and the de Rham cohomology. The
construction leading to short exact exponential sheaf sequences will be presented elsewhere. \\
\\
And now a comment on {\it commutativity}. All of the mentioned extensions of classical smooth differential
geometry, including those in Mallios [1-8], Mallios \& Rosinger, as well as in this paper, have so far been done for
commutative algebras. \\
However, since the recent work of Connes, for instance, much interest has focused on noncommutative structures
as well. It is important, nevertheless, to point out here three facts in this respect. \\
\\
First, when it comes to dealing with singularities, and one does so in a differential context, the approach in
Connes falls far short from reaching the power even of the much earlier linear distribution theory of Schwartz.
Indeed, the only differential type operation in Connes, see [pp 19-28, 287-291], is defined as the commutator with
a fixed operator, that is, a Lie type derivation. In this way, it is a rather particular derivation even within a Banach
algebra. Not to mention that it cannot come anywhere near to dealing with arbitrary closed nowhere dense singularities,
let alone, arbitrary dense singularities which are only restricted by having their complementary dense as well. \\
\\
Second, the existence of noncommutative theories need not at all mean the loss of interest in, let alone, the
abandonment  of commutative theories. Indeed, in many important problems the latter turn out to be both more
effective and far more simple. Not to mention that, so often, noncommutative theories are such only on their so
called 'global' level, while in the last instance of their detailed computations, that is, 'locally', they get reduced to
the commutative. Such a reduction is precisely the reason why they may become tractable, as ways are found to
have their noncommutative complications reduced to commutative computations, as illustrated quite clearly by
matrix theory, for instance. In this connection, one might also refer to N. Bohr's own words, see {\it Bohr's
Correspondence Principle}, according to which 'the description of our measurements of a quantum system must use
classical, commutative ${\bf C}$-numbers.' \\
The present paper, as well, can be seen as another illustration of the usefulness, and also relative simplicity -
based on a setup of sheaves of algebras - of commutative theories, when it comes to the treatment of by far
the largest sets of singularities so far in the literature. \\
\\
Finally, it should be noted that the differential algebras of generalized functions in Rosinger [1-18], Colombeau,
and thus those in Mallios \& Rosinger, or also in this paper, can naturally have their {\it noncommutative}
versions as well. Indeed, such versions are obtained as soon as the respective algebras are constructed not
starting with real or complex valued functions, but with functions which take values in appropriate noncommutative
topological algebras. \\

\bigskip
{\bf 1.2. The Differential Triad $(X, {\cal A}, \partial)$}

\bigskip
It will be useful to recall some basics of the abstract differential
geometry introduced in Mallios [1]. Its initial structure is given by an arbitrary topological space $X$, an
associative, commutative and unital sheaf ${\cal A}$ of ${\bf R}$-algebras on $X$, which in the abstract theory is
the {\it structure sheaf of coefficients}, and finally, by a mapping

\bigskip
(1.1) \quad $ \partial : {\cal A} \longrightarrow \Omega^1 $

\medskip
which is the analog of the usual differential operator, and in the general case is only supposed to have the
following properties : $\Omega^1$ is any sheaf of ${\cal A}$ modules on $X$, while $\partial$ is an ${\bf R}$-linear
sheaf morphism which satisfies the Leibniz rule of product derivative, namely, for any open $U \subseteq X$ and
$\alpha, \beta \in {\cal A}(U), \lambda, \mu \in {\bf R}$, we have

\bigskip
(1.2) \quad $ \partial(\lambda \alpha + \mu \beta) = \lambda \partial \alpha + \mu \partial \beta $

\bigskip
(1.3) \quad $ \partial (\alpha \beta) = (\partial \alpha) \beta + \alpha (\partial \beta) $

\medskip
where we note that ${\bf R} \subseteq {\cal A}$, since ${\cal A}$ is unital.

\medskip
We call $(X, {\cal A}, \partial)$ the {\it differential triad} and recall that in Mallios [1] the notation
$({\cal A}, \partial, \Omega)$ was used instead.

\medskip
We also note that the classical differential geometric setup follows when $X$ is an open subset of ${\bf R}^n$, or
it is an $n$-dimensional manifold, while ${\cal A} = {\cal C}^\infty(X)$ and $\partial$ is the corresponding usual
differential, with $\Omega^1$ being the set of 1-forms.

\medskip
Given now a differential triad $(X, {\cal A}, \partial)$, then for each $n \in {\bf N}, n \geq 2$, we can construct the
$n$-fold exterior product

\bigskip
(1.4) \quad $ \Omega^n = \bigwedge^n \Omega^1 $

\medskip
where the exterior product $\bigwedge$, and its iterates, are constructed with respect to the underlying sheaf of
algebras ${\cal A}$.

\medskip
At this stage, we need to introduce one further entity, namely, we assume the existence of an ${\bf R}$-linear
sheaf morphism

\bigskip
(1.5) \quad $ d^1 : \Omega^1 \longrightarrow \Omega^2 $

\medskip
which satisfies the respective version of the Leibniz rule of product derivative, namely

\bigskip
(1.6) \quad $ d^1(\alpha s) = \alpha (d^1 s) + (\partial \alpha) \bigwedge s $

\medskip
for every $\alpha \in {\cal A}(U), s \in \Omega^1(U)$, and open $U \subseteq X$. Also, we require that

\bigskip
(1.7) \quad $ d^1 \circ d^0 = 0 $

\medskip
where for convenience we denote

\bigskip
(1.8) \quad $ d^0 = \partial $

\medskip
Based on the above, we can now construct the ${\bf R}$-linear sheaf morphism

\bigskip
(1.9) \quad $ d^2 : \Omega^2 \longrightarrow \Omega^3 $

\medskip
by defining it according to

\bigskip
(1.10) \quad $ d^2(s \bigwedge t) = s \bigwedge (d^1 t) + (d^1 s) \bigwedge t $

\medskip
for $s, t \in \Omega^1(U)$ and any open $U \subseteq X$.

\medskip
Finally, as a last assumption, we require that $d^2$ satisfy

\bigskip
(1.11) \quad $ d^2 \circ d^1 = 0 $

\medskip
Based on the above, we can now construct all the ${\bf R}$-linear sheaf morphisms

\bigskip
(1.12) \quad $ d^n : \Omega^n \longrightarrow \Omega^{n+1},~~~ n \in {\bf N},~~ n \geq 3 $

\medskip
by defining them according to

\bigskip
(1.13) \quad $ d^n(s \bigwedge t) = (d^{n-1}s) \bigwedge t + (-1)^{n-1}s \bigwedge (d^1 t) $

\medskip
where $s \in \Omega^{n-1}(U), t \in \Omega^1(U)$, and $U \subseteq X$ is open. \\

\bigskip
{\bf 1.3. De Rham Complexes}

\bigskip
An important fact, which allows the construction of de Rham complexes in the above
abstract framework, is that, within the mentioned constructions, we obtain, see Mallios [1, chap. viii, sect. 8]

\bigskip
{\bf Lemma 1}

\medskip
The relations hold

\bigskip
(1.14) \quad $ d^3 \circ d^2 = d^4 \circ d^3 = ~.~.~.~ = d^{n+1} \circ d^n = ~.~.~.~ = 0,~~~ n \in {\bf N},~~ n \geq 2 $ \\ \\

\bigskip
{\bf 2.~ Sheaves of Algebras of Generalized Functions with Dense Singularities, or Space-Time Foam Algebras} \\

\bigskip
{\bf 2.1. Families of Dense Singularities in Euclidean Spaces}.

\bigskip
Let our domain for generalized functions be any nonvoid open subset $X$ of ${\bf R}^n$. We shall consider various
families of singularities in $X$, each such family being given by a corresponding set ${\cal S}$ of subsets
$\Sigma \subset  X$, with each such subset $\Sigma$ describing a possible set of singularities of a certain given
generalized function. \\
The {\it largest} family of singularities $\Sigma \subset X$ which we can deal with is given by

\bigskip
(2.1) \quad $ {\cal S}_{\cal D}(X) ~=~ \{~ \Sigma \subset X ~~|~~ X \setminus \Sigma ~~\mbox{is dense in}~~ X ~\} $

\medskip
The various families ${\cal S}$ of singularities $\Sigma \subset X$ which we shall use, will therefore each satisfy the
condition ${\cal S} \subseteq {\cal S}_{\cal D}(X)$, see for details Rosinger [14-18]. Among other ones, two such
families which will be of interest are the following

\bigskip
(2.2) \quad $ {\cal S}_{nd}(X) ~=~ \{~ \Sigma \subset X ~~|~~ \Sigma ~~\mbox{is closed and nowhere dense in}~~ X ~\} $

\medskip
and

\bigskip
(2.3) \quad $ {\cal S}_{Baire~I}(X) ~=~ \{~ \Sigma \subset X ~~|~~ \Sigma ~~\mbox{is of first Baire category in}~~ X ~\} $

\medskip
Obviously

\bigskip
(2.4) \quad $ {\cal S}_{nd}(X) ~\subset~ {\cal S}_{Baire~I}(X) ~\subset~ {\cal S}_{\cal D}(X) $

\medskip
In Mallios \& Rosinger, only singularities $\Sigma$ in ${\cal S}_{nd}(X)$ were considered. Thus, in view of (2.4) alone,
it is clear how much more powerful the corresponding results in this paper are, since now we can consider all the
singularities in ${\cal S}_{Baire~I}(X)$, and in fact, even in ${\cal S}_{\cal D}(X)$. \\

\bigskip
{\bf 2.2. Asymptotically Vanishing Ideals}

\bigskip
The construction of the space-time foam algebras, first introduced in Rosinger [14-18], has {\it two} basic ingredients
involved. First, we take any family ${\cal S}$ of singularity sets $\Sigma \subset X$, family which  satisfies the
conditions

\bigskip
(2.5) \quad $ \begin{array}{l}
                                  \forall~ \Sigma \in {\cal S} ~: \\
                                  \\
                                  ~~~ X \setminus \Sigma ~~\mbox{is dense in}~~ X
                      \end{array} $

\medskip
and

\bigskip
(2.6)  \quad $ \begin{array}{l}
                                   \forall~ \Sigma,~ \Sigma^\prime \in {\cal S} ~: \\
                                   \\
                                   \exists~ \Sigma^{\prime \prime} \in {\cal S} ~: \\
                                   \\
                                   ~~~ \Sigma \cup \Sigma^\prime \subseteq \Sigma^{\prime \prime}
                      \end{array} $

\medskip
Clearly, we shall have the inclusion ${\cal S} \subseteq {\cal S}_{\cal D}(X)$ for any such family ${\cal S}$. Also, it is
easy to see that both families ${\cal S}_{nd}(X)$ and ${\cal S}_{Baire~I}(X)$ satisfy conditions (2.5) and (2.6). \\
Now, as the second ingredient, and so far independently of ${\cal S}$ above, we take any right directed partial order
$L = ( \Lambda, \leq )$. In other words, L is such that for each $\lambda,~ \lambda^\prime \in \Lambda$, there exists
$\lambda^{\prime \prime} \in \Lambda$, for which we have $\lambda,~\lambda^\prime \leq \lambda^{\prime \prime}$. The
choice of $L$ may at first appear to be independent of ${\cal S}$, yet in certain specific instances the two may be
related, with the effect that $\Lambda$ may have to be large, see for details subsection 2.8. and Rosinger [15,
subsection 2.6.]. \\
\\
Although we shall only be interested in singularity sets $\Sigma \in {\cal S}_{\cal D}(X)$, the following ideal can be
defined for any $\Sigma \subseteq X$. Indeed, let us denote by

\bigskip
(2.7) \quad $ {\cal J}_{L, \Sigma}(X) $

\medskip
the {\it ideal} in $({\cal C}^\infty(X))^\Lambda$ of all the sequences of smooth functions indexed by $\lambda \in
\Lambda$, namely, $w = (~ w_\lambda ~|~ \lambda \in \Lambda ~) \in ({\cal C}^\infty(X))^\Lambda$, sequences which {\it
outside} of the singularity set $\Sigma$ will satisfy the {\it asymptotic vanishing} condition

\bigskip
(2.8) \quad $ \begin{array}{l}
                                   \forall~ x \in X \setminus \Sigma ~: \\
                                   \\
                                   \exists~ \lambda \in \Lambda ~: \\
                                   \\
                                   \forall~ \mu \in \Lambda,~ \mu \geq \lambda ~: \\
                                   \\
                                   \forall~ p \in {\bf N}^n ~: \\
                                   \\
                                    ~~~D^p w_\mu(x) = 0
                      \end{array} $

\medskip
This means that the sequences of smooth functions $w = ( w_\lambda ~|~ \lambda \in \Lambda )$ in the ideal ${\cal J}_
{L, \Sigma}(X)$ will in a way {\it cover} with their support the singularity set $\Sigma$, and at the same time, they
vanish outside of it, together with all their partial derivatives. \\
In this way, the ideal ${\cal J}_{L, \Sigma}(X)$ carries in an {\it algebraic} manner the information on the singularity set
$\Sigma$. Therefore, a {\it quotient} in which the factorization is made with such ideals may  in certain ways {\it do
away with} singularities, and do so through purely algebraic means, see (2.11), (5.1) below. \\
\\
We note that the assumption about $L = (\Lambda, \leq)$ being right directed is used in proving that ${\cal J}_{L,
\Sigma}(X)$ is indeed an ideal, more precisely that, for $w,~ w^\prime \in {\cal J}_{L, \Sigma}(X)$, we have $w + w^\prime
\in {\cal J}_{L, \Sigma}(X)$. \\
Now, it is easy to see that for $\Sigma,~ \Sigma^\prime \subseteq X$, we have

\bigskip
(2.9) \quad $ \Sigma \subseteq \Sigma^\prime ~\Longrightarrow~ {\cal J}_{L, \Sigma}(X) \subseteq
                                                                                                                            {\cal J}_{L, \Sigma^\prime}(X) $

\medskip
in this way, in view of (2.6), it follows that

\bigskip
(2.10) \quad $ {\cal J}_{L, {\cal S}}(X) ~=~ \bigcup_{\Sigma \in {\cal S}}~ {\cal J}_{L, \Sigma}(X) $

\medskip
is also an {\it ideal} in $({\cal C}^\infty(X))^\Lambda$. \\
\\
It is important to note that for suitable choices of the right directed partial orders $L$, the ideals ${\cal J}_{L,
\Sigma}(X)$, with $\Sigma \in {\cal S}_{\cal D}(X)$, are {\it nontrivial}, that is, they do not collapse to $\{~ 0 ~\}$. Thus in
view of (2.10), the same will hold for the ideals ${\cal J}_{L, {\cal S}}(X)$. In Rosinger [15, section 2] further details are
presented in the case of general singularity  sets $\Sigma \in {\cal S}_{\cal D}(X)$, when the respective right directed
partial orders $L$ which give the nontriviality of the ideals ${\cal J}_{L, \Sigma}(X)$ prove to be rather large, and in
particular, uncountable. In subsection 2.8., we shall show that in the case of the singularity sets $\Sigma$ in
${\cal S}_{nd}(X)$, and in fact, even in ${\cal S}_{\delta~Baire~I}(X)$ - which, see (2.28) below, is a suitable subset of
${\cal S}_{Baire~I}(X)$ - one can have the ideals ${\cal J}_{L, \Sigma}(X)$ nontrivial even for $L = {\bf N}$, respectively,
 for $L = {\bf N} \times {\bf N}$, that is, with $L$ still countable. \\
On the other hand, in view of (2.13), (2,15) below, the mentioned ideals cannot become too large either. \\

\bigskip
{\bf 2.3. Foam Algebras}

\bigskip
In view of the above, for $\Sigma \subseteq X$, we can define the algebra

\bigskip
(2.11) \quad $ B_{L, \Sigma}(X) ~=~ ({\cal C}^\infty(X))^\Lambda / {\cal J}_{L, \Sigma}(X) $

\medskip
However, we shall only be interested in singularity sets $\Sigma \in {\cal S}_{\cal D}(X)$, that is, for which
$X \setminus \Sigma$ is dense in $X$, and in such a case the corresponding algebra $B_{L, \Sigma}(X)$ will be called
a {\it foam algebra}. \\

\bigskip
{\bf 2.4. Multi-Foam Algebras}

\bigskip
With the given family ${\cal S}$ of singularities, based on (2.10), we can associate the {\it multi-foam algebra}

\bigskip
(2.12) \quad $ B_{L, {\cal S}}(X) ~=~ ({\cal C}^\infty(X))^\Lambda / {\cal J}_{L, {\cal S}}(X) $ \\

\bigskip
{\bf 2.5. Space-Time Foam Algebras}

\bigskip
The foam algebras and the multi-foam algebras introduced above will for the sake of simplicity be called together
{\it space-time foam algebras}. Clearly, if the family ${\cal S}$ of singularities consists of one single singularity set
$\Sigma \in {\cal S}_{\cal D}(X)$, that is, ${\cal S} = \{~ \Sigma ~\}$, then conditions (2.5), (2.6) are satisfied, and in this
particular case the concepts of foam and multi-foam algebras are identical, in other words, $B_{L, \{~ \Sigma ~\}}(X) =
B_{L, \Sigma}(X)$. This means that the concept of multi-foam algebra is more general than that of a foam algebra. \\
\\
It is clear from their quotient construction that the space-time foam algebras are associative and commutative.
However, the above constructions can easily be extended to the case when, instead of real valued smooth functions,
we use smooth functions with values in a suitable {\it topological algebra}, or even in an arbitrary {\it normed algebra}.
In such a case the resulting space-time foam algebras will still be associative, but in general they may be
noncommutative, depending on the algebras chosen for the ranges of the smooth functions. \\

\bigskip
{\bf 2.6. Space-Time Foam Algebras as Algebras of Generalized Functions}

\bigskip
The reason why we restrict ourself to singularity sets $\Sigma \in {\cal S}_{\cal D}(X)$, that is, to subsets $\Sigma
\subset X$ for which $X \setminus \Sigma$ is dense in $X$, is due to the implication

\bigskip
(2.13) \quad $ X \setminus \Sigma ~~\mbox{is dense in}~~ X ~~\Longrightarrow~~ {\cal J}_{L, \Sigma}(X) \cap
                                                                                                                              {\cal U}^\infty_\Lambda(X) ~=~ \{~ 0 ~\} $

\medskip
where ${\cal U}^\infty_\Lambda(X)$ denotes the {\it diagonal} of the power $({\cal C}^\infty(X))^\Lambda$, namely, it is
the set of all $u(\psi) = ( \psi_\lambda ~|~ \lambda \in \Lambda )$, where $\psi_\lambda = \psi$, for $\lambda \in
\Lambda$, while $\psi$ ranges over ${\cal C}^\infty(X)$. In this way, we have the algebra isomorphism
${\cal C}^\infty(X) \ni \psi \longmapsto u(\psi) \in {\cal U}^\infty_\Lambda(X)$. \\
\\
The implication (2.13) results immediately from the asymptotic vanishing condition (2.8), and it means that the ideal
${\cal J}_{L, \Sigma}(X)$ is {\it off diagonal}. \\
Yet the importance of (2.13) is that, for $\Sigma \in {\cal S}_{\cal D}(X)$, it gives the following {\it algebra embedding} of
the smooth functions into foam algebras

\bigskip
(2.14) \quad $ {\cal C}^\infty(X) \ni \psi \longmapsto u(\psi) + {\cal J}_{L, \Sigma}(X) \in B_{L, \Sigma}(X) $

\medskip
Now in view of (2.10), it is easy to see that (2.13) will yield the {\it off diagonality} property as well

\bigskip
(2.15) \quad $ {\cal J}_{L, {\cal S}}(X) \cap {\cal U}^\infty_\Lambda(X) ~=~ \{~ 0 ~\} $

\medskip
and thus similarly with (2.14), we obtain the {\it algebra embedding} of smooth functions into multi-foam algebras

\bigskip
(2.16) \quad $ {\cal C}^\infty(X) \ni \psi \longmapsto u(\psi) + {\cal J}_{L, {\cal S}}(X) \in B_{L, {\cal S}}(X) $

\medskip
The algebra embeddings (2.14), (2.16) mean that the foam and multi-foam algebras are in fact {\it algebras of
generalized functions}. Also they mean that the multi-foam algebras are unital, with the respective unit elements
$u(1) + {\cal J}_{L, \Sigma}(X),~ u(1) + {\cal J}_{L, {\cal S}}(X)$. \\
\\
Further, the asymptotic vanishing condition (2.8) also implies quite obviously that, for $\Sigma \subseteq X$, we have

\bigskip
(2.17) \quad $ D^p~ {\cal J}_{L, \Sigma}(X) \subseteq  {\cal J}_{L, \Sigma}(X), ~~~\mbox{for}~~ p \in {\bf N}^n $

\medskip
where $D^p$ denotes the termwise $p$-th order partial derivation of sequences of smooth functions, applied to
each such sequence in the ideal ${\cal J}_{L, \Sigma}(X)$. \\
Then again, in view of (2.10), we obtain

\bigskip
(2.18) \quad $ D^p~ {\cal J}_{L, {\cal S}}(X) \subseteq  {\cal J}_{L, {\cal S}}(X), ~~~\mbox{for}~~ p \in {\bf N}^n $

\medskip
Now (2.17), (2.18) mean that the the foam and multi-foam algebras are in fact {\it differential algebras}, namely

\bigskip
(2.19) \quad $ D^p~ B_{L, \Sigma}(X) \subseteq B_{L, \Sigma}(X), ~~~\mbox{for}~~ p \in {\bf N}^n $

\medskip
where $\Sigma \in {\cal S}_{\cal D}(X)$, while we also have

\bigskip
(2.20) \quad $  D^p~ B_{L, {\cal S}}(X) \subseteq B_{L, {\cal S}}(X), ~~~\mbox{for}~~ p \in {\bf N}^n $

\medskip
In this way we obtain that the foam and multi-foam algebras are {\it differential algebras of generalized functions}. \\
Also, the multi-foam algebras contain the Schwartz distributions, that is, we have the {\it linear embeddings} which
respect the arbitrary partial derivation of smooth functions

\medskip
(2.21) \quad $ {\cal D}^\prime(X) \subset B_{L, \Sigma}(X),~~~ \mbox{for}~~ \Sigma \in {\cal S}_{\cal D}(X) $

\bigskip
(2.22) \quad $  {\cal D}^\prime(X) \subset B_{L, {\cal S}}(X) $

\medskip
Indeed, we can recall the wide ranging purely algebraic characterization of all those quotient type algebras of
generalized functions in which one can embed linearly the Schwartz distributions, a characterization first given in
1980, see Rosinger [4, pp. 75-88], Rosinger [5, pp. 306-315], Rosinger [6, pp. 234-244]. According to that
characterization - which also contains the Colombeau algebras as a particular case - the necessary and sufficient
condition for the existence of the linear embedding (2.21) is precisely the off diagonality condition in (2.13). Similarly,
the necessary and sufficient condition for the existence of the linear embedding (2.22) is exactly the off diagonality
condition (2.15). \\
\\
One more property of the foam and multi-foam algebras will prove to be useful. Namely, in view of (2.10) it is clear that,
for every $\Sigma \in {\cal S}$, we have the inclusion ${\cal J}_{L, \Sigma}(X) \subseteq {\cal J}_{L, {\cal S}}$, and thus
we obtain the {\it surjective algebra homomorphism}

\bigskip
(2.23) \quad $ B_{L, \Sigma}(X) \ni w + {\cal J}_{L, \Sigma}(X) \longmapsto w + {\cal J}_{L, {\cal S}}(X) \in
                                                                                                                                               B_{L, {\cal S}}(X) $

\medskip
As we shall see in the next subsection, (2.23) can naturally be interpreted as meaning that the typical generalized
functions in $B_{L, {\cal S}}(X)$ are {\it more regular} than those in $B_{L, \Sigma}(X)$. \\

\bigskip

{\bf 2.7. Regularity of Generalized Functions}

\medskip
One natural way to interpret (2.23) in the given context of generalized functions is the following. \\
\\
Given two spaces of generalized functions $E$ and $F$, such as for instance

\bigskip
(2.24) \quad $ {\cal C}^\infty(X) \subset E \subset F $

\medskip
then the {\it larger} the space $F$, the {\it less regular} its typical elements can appear to be, when compared with
those of $E$. By the same token, the {\it smaller} the space $E$, the {\it more regular}, compared to those of $F$, one
can consider its typical elements. \\
\\
This kind of regularity we can call {\it subset regularity}. \\
\\
On the other hand, given a {\it surjective} mapping

$$ E \rightarrow F $$

\medskip
it may at first sight appeaar that one could again consider that the typical elements of $F$ are at least as {\it regular} as
those of $E$. \\
However, as the following simple example shows it, additional arguments may be needed for such a conclusion.
Indeed, we clearly have, for instance, the inclusions

$$ {\cal C}^\infty ( {\bf R} ) ~\subset~ {\cal C}^1 ( {\bf R} ) ~\subset~ {\cal C}^0 ( {\bf R} ) $$

\medskip
as well as the surjective linear mapping given by the usual derivative, namely

$$ D : {\cal C}^1 ( {\bf R} ) ~\longrightarrow~ {\cal C}^0 ( {\bf R} ) $$

\medskip
yet it is the elements of ${\cal C}^1 ( {\bf R} )$ which are considered to be more regular than those of ${\cal C}^0
( {\bf R} )$. \\
In this way, in order to be able to support the argument that in the case of a surjective mapping $E \longrightarrow F$,
we can indeed say about $F$ to have more regular elements than those of $E$, the respective surjective mapping
should enjoy certain suitable additional properties. And clearly, such is not the case with the derivative mapping in the
counterexample above. \\
\\
However, if both $E$ and $F$ are {\it quotient vector spaces} of the form specified next, and the surjective mapping is
the {\it canonical} one between them, namely

\bigskip
(2.25) \quad $ E ~=~ {\cal S} / {\cal V} \ni s ~+~ {\cal V} ~\longmapsto~ s ~+~ {\cal W} \in F ~=~ {\cal S} / {\cal W} $

\medskip
where ${\cal V} \subseteq {\cal W} \subseteq {\cal S}$ are vector spaces, then one can see the elements of $F$ as
being more {\it regular} then those of $E$, since ${\cal W}$ may factor out in $F$ more than does ${\cal V}$ in $E$. \\
\\
This kind of regularity we shall call {\it quotient regularity}. \\
\\
In this way, view of (2.23), we can consider that, owing to the given {\it surjective} algebra homomorphism, the typical
elements of the multi-foam algebra $B_{L, {\cal S}}(X)$ can be seen as being more {\it quotient regular} than the typical
elements of the foam algebra $B_{L, \Sigma}(X)$. \\
Indeed, the algebra $B_{L, {\cal S}}(X) $ is obtained by factoring the same $({\cal C}^\infty(X))^\Lambda$ as in the case
of the algebra $B_{L, \Sigma}(X)$, this time however by the significantly {\it larger} ideal ${\cal J}_{L, {\cal S}}(X)$, an
ideal which, unlike any of the individual ideals ${\cal J}_{L, \Sigma}(X)$, can simultaneously deal with {\it all} the
singularity sets $\Sigma \in {\cal S}$, some, or in fact, many of which can be {\it dense} in $X$. Further details related to
the connection between {\it regularization} in the above sense, and on the other hand, properties of {\it stability}, {\it
generality} and {\it exactness} of generalized functions and solutions can be found in Rosinger [4-6]. \\

\bigskip
{\bf 2.8. Nontriviality of Ideals}

\medskip
Let us take any nonvoid singularity set $\Sigma \in {\cal S}_{nd}(X)$. Since $\Sigma$ is closed, we can take a
sequence of nonvoid open subsets $Y_l \subset X$, with $l \in {\bf N}$, such that $\Sigma = \cap_{l \in {\bf N}}~Y_l$.
We can also assume that the $Y_l$ are decreasing in $l$, since we can replace every $Y_l$ with the finite intersection
$\cap_{k \leq l}~ Y_k$. But for each $Y_l$, Kahn, there exists $\alpha_l \in {\cal C}^\infty(X)$, such that $\alpha_l
= 1$ on $\Sigma$, and $\alpha_l = 0$ on $X \setminus Y_l$. Now in view of (2.8) it is easy to check that the resulting
sequence of smooth functions on $X$ satisfies

\bigskip
(2.26) \quad $ \alpha = ( \alpha_l ~|~ l \in {\bf N} ) \in {\cal J}_{{\bf N}, \Sigma}(X) $

\medskip
and clearly, $\alpha$ in not a trivial sequence, since $\phi \neq \Sigma \subseteq ~\mbox{supp}~ \alpha_l$, for $l \in
{\bf N}$. \\
\\
We can note, however, that the above argument leading to (2.26) need not necessarily apply to subsets $\Sigma
\subset X$ which are not closed, but whose closure is nevertheless nowhere dense in $X$. In such a case one can use
the more general method in Rosinger [15, section 2], which will give nontrivial sequences similar to $\alpha$ above,
however, their index sets will no longer be countable.

\medskip
Let us take now any nonvoid singularity set $\Sigma \in {\cal S}_{Baire~I}(X)$. Then there exists a sequence of closed
and nowhere dense subsets $\Sigma_l \subset X$, with $l \in {\bf N}$, such that

\bigskip
(2.27) \quad $ \Sigma ~\subseteq~ \bigcup_{l \in {\bf N}}~ \Sigma_l $

\medskip
where the equality need not necessarily hold. Therefore, let us consider the subset of ${\cal S}_{Baire~I}(X)$ denoted by

\bigskip
(2.28) \quad $ {\cal S}_{\delta~Baire~I}(X) $

\medskip
whose elements are all those singularity sets $\Sigma$ for which we have equality in (2.27). Obviously, we can
assume that the $\Sigma_l$ are increasing in $l$, since we can replace each $\Sigma_l$ with the finite union $\cup_{k
\leq l}~ \Sigma_k$. \\
\\
Given now a nonvoid $\Sigma \in {\cal S}_{\delta~Baire~I}(X)$ and a corresponding representation $\Sigma = \cup_{l \in
{\bf N}}~\Sigma_l$, with suitable closed and nowehere dense subsets $\Sigma_l \subseteq X$ which are increasing in
$l$, we can, as above, find for each $\Sigma_l$ a representation $\Sigma_l = \cap_{k \in {\bf N}}~Y_{l k}$, with nonvoid
open subsets $Y_{l k} \subset X$. Further, we can assume that for $l,~ l^\prime,~ k,~ k^\prime \in {\bf N},~ l \leq l^\prime,~
k \leq k^\prime$, we have $Y_{l k} \supseteq Y_{l^\prime k^\prime}$, since we can replace every $Y_{l k}$ with the finite
intersection $\cap_{l^\prime \leq l, k^\prime \leq k}~ Y_{l^\prime k^\prime}$. Now, for each $Y_{l k}$, we can find
$\alpha_{l k } \in {\cal C}^\infty(X)$, such that $\alpha_{l k} = 1$ on $\Sigma_l$, while $\alpha_{l k} = 0$ on  $X \setminus
Y_{l k}$. \\
Let us now take $L = ( \Lambda, \leq )$, where $\Lambda = {\bf N} \times {\bf N}$ and for $(l,k),~ (l^\prime,k^\prime) \in
\Lambda = {\bf N} \times {\bf N}$ we set $(l,k) \leq (l^\prime,k^\prime)$, if and only if $l \leq l^\prime$ and $k \leq
k^\prime$. Then (2.8) will easily give

\bigskip
(2.29) \quad $ \alpha = ( \alpha_{l k} ~|~ (l,k) \in {\bf N} \times {\bf N} ) \in {\cal J}_{{\bf N} \times {\bf N},
                                                                                                        \Sigma}(X) $

\medskip
And again, $\alpha$ is not a trivial sequence, since $\phi \neq \Sigma \subseteq \cup_{l \in {\bf N}}~ \mbox{supp}~
\alpha_{l, k_l}$, for every given choice of $k_l \in {\bf N}$, with $l \in {\bf N}$. \\
\\
In case our singularity set $\Sigma$ belongs to ${\cal S}_{Baire~I}(X)$ but not to ${\cal S}_{\delta~Baire~I}(X)$,
then the above approach need no longer work. However, we can still apply the mentioned more general method in
Rosinger [15, section 2], in order to construct nontrivial sequences in ${\cal J}_{L, \Sigma}(X)$ although this time the
corresponding index sets $\Lambda$ may be uncountable. \\

\bigskip
{\bf 2.9. Relations between Algebras with the Same Singularities}

\bigskip
The above, and especially subsection 2.8., leads to the following question. Let us assume given a certain nonvoid
singularity set $\Sigma \in {\cal S}_{\cal D}(X)$. If we now consider two right directed partial orders $L = (\Lambda,
\leq)$ and $L^\prime = (\Lambda^\prime, \leq)$, is there then any relevant relationship between the corresponding two
foam algebras

\bigskip
(2.30) \quad $ B_{L, \Sigma}(X) ~~~\mbox{and}~~~ B_{L^\prime, \Sigma}(X) ~~ ? $

\medskip
A rather simple positive answer can be given in the following particular case. Let us assume that $\Lambda$ is a {\it
cofinal} subset of $\Lambda^\prime$, that is, the partial order on $\Lambda$ is induced by that on $\Lambda^\prime$,
and in addition, we also have satisfied the condition

\bigskip
(2.31) \quad $ \begin{array}{l}
                             \forall~  \lambda^\prime \in \Lambda^\prime ~: \\
                             \\
                             \exists~ \lambda \in \Lambda ~: \\
                             \\
                             ~~~ \lambda^\prime \leq \lambda
                       \end{array} $

\medskip
Then considering the surjective algebra homomorphism

\bigskip
(2.32) \quad $ \begin{array}{l}
                      ({\cal C}^\infty(X))^{\Lambda^\prime} \ni s^\prime = ( s^\prime_{\lambda^\prime} ~|~ \lambda^\prime \in
                      \Lambda^\prime ) ~\stackrel{\rho}{\longmapsto}~ \\
                      \\
                       ~~~~~~~~~~~~~~~~~~~~~~~~~~~~~~~~~~~~~~~~~~s = ( s^\prime_{\lambda^\prime} ~|~ \lambda^\prime
                                                                                                              \in \Lambda ) \in ({\cal C}^\infty(X))^{\Lambda}
                      \end{array} $

\medskip
and based on (2.8), one can easily note the property

\bigskip
(2.33) \quad $ \rho~ {\cal J}_{\Lambda^\prime, \Sigma}(X) ~\subseteq~  {\cal J}_{\Lambda, \Sigma}(X) $

\medskip
In this way, one can obtain the {\it surjective} algebra homomorphism of foam algebras, given by

\bigskip
(2.34) \quad $ B_{\Lambda^\prime, \Sigma}(X) \ni s^\prime + {\cal J}_{\Lambda^\prime, \Sigma}(X) ~\stackrel{\rho}
                               {\longmapsto}~ \rho~ s^\prime + {\cal J}_{\Lambda, \Sigma}(X) \in B_{\Lambda, \Sigma}(X) $

\medskip
In the terms of the interpretation in subsection 2.7., the meaning of (2.34) is that the foam algebra $B_{\Lambda,
\Sigma}(X)$ has its typical generalized functions {\it more} regular than those of  $B_{\Lambda^\prime, \Sigma}(X)$.
Thus in such terms, foam algebras which correspond to a {\it smaller cofinal} partial order $L$, can be seen as {\it
more} regular. \\
However, there may be many other kind of relationships between two partial orders $L$ and $L^\prime$, such as for
instace in the case of ${\bf N}$ in (2.26), and ${\bf N} \times {\bf N}$ in (2.29). Therefore the problem in (2.30) may in
general present certain difficulties. \\
\\
Needless to say, similar results and comments hold in the case of the space-time foam algebras. \\

\bigskip
{\bf 2.10. The Flabby and Fine Sheaf Property}

\bigskip
We recall that in the abstract differential geometry in Mallios [1], the structure coefficients - no longer given by smooth
functions - are sheaves of algebras. And for cohomological and then differential reasons, it proves to be very profitable
for such sheaves ( loc.cit. ) to  be {\it flabby} and/or {\it fine}. \\
With the space-time foam algebras of generalized functions presented in subsections 2.1. - 2.9., the issue of being fine
sheaves should, in principle, not raise difficulties, since these algebras are constructed by using classes of
sequences of smooth functions. Also we can recall that many other algebras or spaces of generalized functions in the
literature prove to be fine sheaves. \\
However, the issue of flabbiness is a priori not so obvious. And it is even less so, if we recall that most of the familiar
spaces of generalized functions - and that includes the Scwartz distributions and the Colombeau algebras, among
others - fail to be flabby sheaves. Moreover, their lack of flabbiness is quite closely related to a number of deficiencies,
as shown for instance in Kaneko. \\
\\
Fortunately, as with the algebras of generalized functions used in Mallios \& Rosinger, which could deal with arbitrary
closed nowehere dense singularities, so with the space-time foam algebras used in this paper, which can deal with the
much large class of dense singularities in ${\cal S}_{\cal D}(X)$, they prove to be flabby sheaves, as well. \\
\\
Let us first define a large class of space-time foam algebras $B_{L, {\cal S}}(X)$ on nonvoid open subsets $X
\subseteq {\bf R}^n$, each of which will, in Lemma 2 next, prove to have a {\it fine sheaf} structure. From the proof it will
also follow that the respective algebras are {\it flabby sheaves} as well, in case they satisfy a further rather natural
condition. This class contains the nowehre dense algebras, and thus the result presented here is a significant
extension of the similar recent result in Mallios \& Rosinger [Lemma 2], which was fundamental for that paper. \\
\\
Given a family ${\cal S}$ of singularity sets $\Sigma \subset X$ for which the conditions (2.5), (2.6) hold, we call that
family {\it locally finitely additive}, if and only if it satisfies also the condition : \\
For any sequence of singularity sets $\Sigma_l \in {\cal S}$, with $l \in {\bf N}$, if we take $\Sigma = \bigcup_{l \in
{\bf N}}~\Sigma_l$, then for every nonvoid open subset $U \subseteq X$, we have $\Sigma \cap U \in {\cal S}|_U$,
 whenever

\bigskip
(2.35) \quad $ \begin{array}{l}
                           \forall~ x \in U ~: \\
                           \\
                           \exists~ \Delta \subseteq U,~ \Delta ~~\mbox{neighbourhood of}~ x ~: \\
                           \\
                           ~~~ \{~ l \in {\bf N} ~|~ \Sigma_l \cap \Delta \neq \phi ~\}~~~ \mbox{is a finite set of indices}
                     \end{array} $

\medskip
It is easy to verify that, see (2.2), (2.3), the families of singularities ${\cal S}_{nd}(X)$ and ${\cal S}_{Baire~I}(X)$ are
both locally finitely additive. \\
Indeed, ${\cal S}_{Baire~I}(X)$ is trivially so, since any countable union of first Baire category sets is still of first Baire
category. As far as ${\cal S}_{nd}(X)$ is concerned, it suffices to note two facts. First, a subset of a topological space is
closed and nowehre dense, if and only if it is {\it locally}so, that is, in the neighbourhood of evey point. Second, a finite
union of closed nowhere dense sets is again closed and nowhere dense. \\
\\
Let us also note that ${\cal S}_{nd}(X)$ and ${\cal S}_{Baire~I}(X)$ are among those classes of singularities ${\cal S}$
which for every nonvoid open subset $U \subseteq X$, satisfy the condition

\bigskip
(2.36) \quad $ {\cal S}|_U ~\subseteq~ {\cal S} $

\medskip
where we defined the restriction ${\cal S}|_U$ of ${\cal S}$ to $U$, according to

\bigskip
(2.37) \quad $ {\cal S}|_U ~=~ \{~ \Sigma \cap U ~|~ \Sigma \in {\cal S} ~\} $

\medskip
We shall use the concept of {\it sheaf} as is defined by its {\it sections}, see Bredon, or Mallios [1],
Oberguggenberger \& Rosinger, Mallios \& Rosinger. In particular, we shall deal with {\it restriction} mappings to
nonvoid open subsets $U \subseteq X$. \\
\\
Let us assume given a family ${\cal S}$ of singularities which satisfies the conditions (2.5), (2.6). Then it is clear
that for every nonvoid open subset $U \subseteq X$, the restriction ${\cal S}|_U$ of ${\cal S}$ to $U$ will also satisfy
(2.5), (2.6), this time on $U$. \\
\\
Let us now be given any right directed partial order $L = ( \Lambda, \leq )$. Then the restriction to nonvoid open
subsets $U \subseteq X$ of the space-time foam algebra $B_{L, {\cal S}}(X)$ is the family of space-time foam algebras

\bigskip
(2.38) \quad $ {\cal B}_{L, {\cal S}, X} ~=~ (~ B_{L, {\cal S}|_U}(U) ~|~ U \subseteq X,~~ U ~\mbox{nonvoid open} ~) $

\medskip
a relation which follows easily, if we take into account (2.37), and the fact that

\bigskip
(2.39) \quad $ B_{L, {\cal S}}(X)|_U = B_{L, {\cal S}|_U}(U) $

\medskip
which is a direct consequence of (2.12), (2.10), as well as of the obvious relation, see (2.8)

\bigskip
(2.40) \quad $ {\cal J}_{L, \Sigma}(X)|_U ~=~ {\cal J}_{L, \Sigma \cap U}(U), ~~~\mbox{for}~~ \Sigma \subseteq X $

\medskip
We can also note that in the case $\Sigma \cap U = \phi$, the ideal ${\cal J}_{L, \phi}(U)$, and thus the algebra
$B_{L, \phi}(U)$ are still well defined, as long as $U$ is open and nonvoid, see (2.8), (2.10), (2.12).

\bigskip
{\bf Lemma 2}

\medskip
Given on a nonvoid open subset $X \subseteq {\bf R}^n$ any family of singularities ${\cal S}$ which is locally
finitely additive. \\
Then for every right directed partial order $L = ( \Lambda, \leq )$, the family of space-time foam algebras, see (2.38)

\bigskip
(2.41) \quad $ {\cal B}_{L, {\cal S}, X} ~=~ (~ B_{L, {\cal S}|_U}(U) ~|~ U \subseteq X,~~ U ~\mbox{nonvoid open} ~) $

\medskip
is a fine sheaf on $X$. \\
\\
If in addition ${\cal S}$ has the property

\bigskip
(2.42) \quad $ \begin{array}{l}
                    \forall~ \Sigma \in {\cal S},~~ U \subseteq X,~~ U ~\mbox{nonvoid open},~~~  \Gamma \in {\cal S}_{nd}(U) ~: \\
                    \\
                    ~~~ ( \Sigma \cap U ) \cup \Gamma \in {\cal S}|_U
               \end{array} $

\medskip
and ${\bf N}$ is, see subsection 2.9., cofinal in $\Lambda$, then ${\cal B}_{L, {\cal S}, X}$ in (2.41) is also a
flabby sheaf on $X$.

\bigskip
{\bf Proof}. See Appendix.

\bigskip
{\bf Note 1}

\medskip
The classes of singularities ${\cal S}_{nd}(X)$ and ${\cal S}_{Baire~I}(X)$ satisfy the conditions in Lemma 2 above .

\bigskip
{\bf Note 2}

\medskip
If we consider ${\bf N} \times {\bf N}$ with the partial order in subsection 2.8., and we embed ${\bf N}$ into
${\bf N} \times {\bf N}$ through the diagonal mapping ${\bf N} \ni l \longmapsto (l,l) \in {\bf N} \times {\bf N}$,
then ${\bf N}$ will be {\it cofinal} in ${\bf N} \times {\bf N}$. Thus in view of the Lemma 2 and Note 1 above, it
follows that

$$ {\cal B}_{{\bf N} \times {\bf N}, {\cal S}_{Baire~I}(X), X} ~=~ (~ B_{L, {\cal S}_{Baire~I}(X)|_U}(U) ~|~ U
                                                                    \subseteq X,~~ U ~\mbox{nonvoid open} ~) $$

\medskip
is a {\it fine} and {\it flabby sheaf}. \\
This result is nontrivial since ${\cal S}_{Baire~I}(X)$ contains lots of singularity sets $\Sigma \subseteq X$, which
are both {\it dense} in $X$ and {\it uncountable}. In particular, this result is a significant strengthening of an earlier
similar result in Mallios \& Rosinger [Lemma 2], where it was only given in the case of the family of singularities
${\cal S}_{nd}(X)$. \\
In Rosinger [16] the above Lemma 2 is in fact proved for $X$ any finite dimensional smooth manifold.

\bigskip
{\bf Note 3}

\medskip
In the context of {\it flabbiness} of spaces of functions or generalized functions, the presence of ${\cal S}_{nd}(X)$ in
condition (2.42) appears to be quite natural. For instance, as seen in Oberguggenberger \& Rosinger [Remark 7.5, pp.
142-146], the class ${\cal S}_{nd}(X)$ of closed nowhere dense singularities appears when one constructs the {\it
smallest flabby sheaf} which contains ${\cal C}^\infty(X)$ for a nonvoid open subset $X \subseteq {\bf R}^n$. The same
happens when constructing the smallest flabby sheaf containing ${\cal C}^0(X)$.

\bigskip
{\bf Note 4}

\medskip
It is useful to note that there are other ways as well than in (2.8) to define ideals which lead to
differential algebras of generalized functions with {\it dense singularities}. Here we mention in short
one such way used recently in Rosinger [18]. Given any family of singularities ${\cal S}$ such as in
(2.5), (2.6), we associate with it the {\it ideal} in $( {\cal C}^\infty ( X ) )^{\bf N}$ defined
by

\bigskip
(2.43) \quad $ {\cal J}^f_{\cal S} ( X ) ~=~ \bigcup_{\Sigma \in {\cal S}}~ {\cal J}^f_\Sigma ( X ) $

\medskip
where ${\cal J}^f_\Sigma ( X )$ is the ideal in $( {\cal C}^\infty ( X ) )^{\bf N}$ of all sequences
of smooth functions $w = ( w_\nu | \nu \in {\bf N} )$, which satisfy the condition

\bigskip
(2.44) \quad $ \begin{array}{l}
                    \forall~~ x \in X \setminus \Sigma,~ l \in {\bf N} ~: \\ \\
                    \exists~~ \nu \in {\bf N} ~: \\ \\
                    \forall~~ \mu \in {\bf N},~ \mu \geq \nu,~ p \in {\bf N}^n,~
                                    | p | \leq l ~: \\ \\
                    ~~~ D^p W_\mu ( x ) ~=~ 0
                \end{array} $

\medskip
Clearly, the {\it asymptotic vanishing} condition in (2.44) is {\it weaker} than that in (2.8).
Yet as seen in Rosinger [18], the resulting ideals, and consequently, differential algebras of
generalized functions can handle a variety of problems related to {\it dense
singularities}. \\ \\

\bigskip
{\bf 3.~ Differential Geometry on Space-Time Foam Algebras of Generalized Functions}

\bigskip
We now show that the {\it abstract} differential geometry in Mallios [1], presented in short in section 1, can be
implemented - as a particular case, which allows {\it dense} singularities - by using as {\it structure sheaf of
coefficients} the sheaf of space-time foam algebras of generalized functions, see (2.41) in Lemma 2, section 2. \\

\bigskip
{\bf 3.1. Space-Time Foam Differential Triads}

\medskip
We construct a large variety of differential triads as follows. We take $X$ any nonvoid open subset of ${\bf R}^n$.
Further, we can choose on $X$ in a variety of ways a family of singularities ${\cal S}$ which
satisfies (2.5), (2.6). Once this is done, then through (2.12), (2.41), we are led to the structure sheaf of coefficients given
by the corresponding sheaf ${\cal B}_{L, {\cal S}, X}$ of space-time foam differential algebras of generalized
functions. \\
At that stage, according to the abstract theory, we have to choose the third element, namely, $\partial$, of the
differential triad, see (1.1). For that purpose, first we define for every nonvoid open $U \subseteq X$ the corresponding
$B_{L, {\cal S}}(U)$-module $\Omega^1(U)$, as being the free $B_{L, {\cal S}} (U)$-module of rank $n$, with free
generators $dx_1, dx_2, dx_3,~.~.~.~, dx_n$. In this way, in view of (2.20), the elements of $\Omega^1(U)$ are given by
all

\bigskip
(3.1) \quad $ \sum_{i=1}^n~ V_i~ dx_i $

\medskip
where $V_i \in B_{L, {\cal S}_L} (U)$. \\
Thus we can define in our specific case the desired third element of the differential triad, namely, the sheaf morphism
$\partial$ in (1.1), and do so according to

\bigskip
(3.2) \quad $ \begin{array}{l}
                          B_{L, {\cal S}} (U) \stackrel{\partial}{\longrightarrow} \Omega^1(U) \\
                          \\
                          ~~~~~~ V \longmapsto \sum_{i=1}^n~ (\partial_i V)~ dx_i
                     \end{array} $

\medskip
where as usual, $\partial_i$ denotes the partial derivation with respect to the $i$-th independent variable $x_i$.

\medskip
The effect of the choice in (3.1), (3.2) is the following easy to prove result

\bigskip
{\bf Lemma 3}

\medskip
$(X,~ {\cal B}_{L, {\cal S}, X},~ \partial )$ is a differential triad.

\hfill $\diamondsuit \diamondsuit \diamondsuit$

\medskip
Now according to the abstract theory in section 1, we are at the stage where we have to define the ${\bf R}$-linear
sheaf morphism $d^1$ in (1.5). In view of (3.1), (3.2) and (1.4), for every nonvoid open $U \subseteq X$, we shall
take

\bigskip
(3.3) \quad $ d^1 : \Omega^1(U) \longrightarrow \Omega^2(U) = \Omega^1(U) \bigwedge \Omega^1(U) $

\medskip
where

\bigskip
(3.4) \quad $ d^1( \sum_{i=1}^n~ V_i~ dx_i ) ~=~ \sum_{j=1}^n~ \sum_{i=1}^n~ ( \partial_j~ V_i )~ dx_j \wedge dx_i $

\medskip
Through a direct computation based on (2.18), (2.20), one can verify that, with this definition in (3.3), (3.4), $d^1$
will indeed satisfy conditions (1.6), (1.7). \\
Finally, by implementing (1.9) through (1.10), and using the fact that it is true for smooth functions,
another direct computation will give (1.11).

\bigskip
{\bf Remark 1}

\medskip
Let us recall that in Rosinger [1-13] a large variety, and in fact, infinitely many classes of differential algebras of
generalized functions were constructed. Also a wide ranging purely algebraic characterization was given there for
those algebras which contain the linear vector space of Schwartz distributions. \\
Untill more recently, only two particular cases of these classes of algebras have been used in the study of global
generalized solutions of linear and nonlinear PDEs. Namely, first was the class of the nowhere dense differential
algebras of generalized functions in Rosinger [3-13], while later came the class of algebras considered in
Colombeau. \\
These latter algebras, since they also contain the Schwartz distributions, are, in view of the above mentioned
algebraic characterization, by necessity a particular case of the classes of algebras of generalized functions first
introduced in Rosinger [1-13]. \\
\\
The Colombeau algebras of generalized functions enjoy a rather simple and direct connection with the Schwartz
distributions, and therefore, with a variety of Sobolev spaces. This led to their relative popularity in the study of
generalized solutions of PDEs. \\
Compared however with the nowhere dense differential algebras of generalized functions, let alone with the space-time
foam differential algebras of generalized functions used in this paper, the Colombeau algebras suffer from several
important limitations. Among them, relevant to this paper, and in general, to abstract differential geometry, is the
following. There are growth conditions which the generalized functions must satisfy in the neighbourhood of their
singularities. The effect, among others, is that the Colombeau algebras - just as the Schwartz distributions, for instance
- do {\it not} form a flabby sheaf. In this way, they would not be the appropriate sheaf of structure coefficients even in
such a general theory as the abstract differential geometry in Mallios [1]. In particular, owing to the growth conditions
they have to satisfy, the Colombeau algebras do not allow exponential short exact sequences to be defined, see
Mallios \& Rosinger. This indeed constitutes a severe shortcoming for important applications in mathematical physics,
for instance, {\it geometric (pre)quantization}, as e.g. Weil's integrality theorem, see Mallios [4], or even Mallios [1,
chap. viii, sect. 11]. \\
\\
On the other hand, the earlier introduced nowhere dense algebras do not suffer from any of the above two
limitations. Indeed, the nowhere dense algebras allow singularities on arbitrary closed nowhere dense sets,
therefore, such singularity sets can have arbitrary large positive Lebesgue measure, Oxtoby. Furhtermore,
in the nowhere dense algebras there are no any conditions asked on generalized functions in the neighbourhood
of their singularities. \\
\\
In this paper, the use of the space-time foam differential algebras of generalized functions, introduced recently in
Rosinger [14-18], brings a further significant enlargement of the possibilities already given by the nowhere dense
algebras, and applied in Mallios \& Rosinger. Indeed, this time the singularities can be concentrated on arbitrary
subsets, including dense ones, provided that their complementary, that is, the set of nonsingular points, is still
dense. Furhtermore, as already in the case of the nowhere dense algebras, also in the space-time foam algebras, no
any kind of condition is asked on the generalized functions in the neighbourhood of their singularities. \\
\\
Finally, it should be noted that, since one of the major interests in differential geometry, including in its abstract
version in Mallios [1-8], comes from general relativity, it is important to have in the respective frameworks strong
and general enough results on the existence of solutions for nonlinear PDEs. In this respect, one could already
obtain in the framework of the nowehre dense algebras a rather general, and in fact, {\it type independent} and
{\it global} version of the classical Cauchy-Kovalevaskaia theorem, see Rosinger [7-9]. Indeed, one can prove
that every analytic nonlinear PDE, with every associated noncharacteristic analytic initial value problem, has a
{\it global} generalized solution, which is analytic on the whole domain of definition of the respective PDE, except
for a closed nowhere dense set, set which can be chosen to have zero Lebesgue measure. \\
This global type independent existence results is, fortunately, preserved in the case of the space-time foam algebras
as well, see Rosinger [14,15]. \\
So far, one could not obtain any kind of similarly general and powerful existence of solutions result in any of the
infinitely many other classes of algebras of generalized functions, including in the Colombeau class. \\ \\

\bigskip
{\bf 4.~ De Rham Cohomology with Dense Singularities}

\bigskip
Let us suppose given a nonvoid open subset $X \subseteq {\bf R}^n$ and any family ${\cal S}$  of singularities
on it, see (2.5), (2.6). \\
The corresponding space-time foam differential triad $(X, {\cal B}_{L, {\cal S}, X}, \partial )$, see Lemma 3,
section 3, leads according to section 1, and the general theory in Mallios [1, chap. iii], to the following complex of
${\bf R}$-linear sheaf morphisms, that is, to the {\it de Rham complex with dense singularities}

\bigskip
(4.1) \quad $ 0 \longrightarrow {\bf R} \stackrel{\epsilon}{\longrightarrow} {\cal B}_{L, {\cal S}, X} \stackrel{\partial}
                                {\longrightarrow} \Omega^1 \stackrel{d^1}{\longrightarrow} \Omega^2 \stackrel{d^2}
                                        {\longrightarrow} ~.~.~.~   $

\medskip
And as an extension of the similar result for nowhere dense singularities in Mallios \& Rosinger, here we have for
dense singularities

\bigskip
{\bf Theorem 1}

\medskip
The de Rham complex (4.1) is exact, namely,

\bigskip
(4.2) \quad $ \mbox{ker}~ d^{n+1} ~=~ \mbox{im}~ d^n,~~~ n \in {\bf N} $

\medskip
where we denoted, see (1.8), $d^0 = \partial$.

\medskip
{\bf Proof}. We note that (4.2) is equivalent with the Poincare Lemma in the abstract differential geometry which
corresponds to ${\cal B}_{L, {\cal S}, X}$, or more precisely, to the space-time foam differential triad
$(X, {\cal B}_{L, {\cal S}, X}, \partial)$, see Mallios [1]. This means that, locally, every closed differential form is exact.
In this way, the exactness of (4.1) can be checked 'fiberwise', that is, through the use of (2.18), (2.20), in other words,
by reducing it to the classical case of smooth functions. \\
As $X$ is a nonvoid open subset of ${\bf R}^n$, it follows that $X$ is Hausdorff and paracompact. And in view of
Lemma 2, section 2, ${\cal B}_{L, {\cal S}, X}$ is a fine sheaf on $X$, thus the same holds for all $\Omega^{n+1}$,
with $n \in {\bf N}$, see Mallios [1, chap. iii, (8.56) ].

\hfill $\diamondsuit \diamondsuit \diamondsuit$

\medskip
Now, in the terms of Mallios [1, chap. iii, (8.24) ], one obtains

\bigskip
{\bf Corollary 1}

\medskip
The complex (4.1) provides a fine, hence, a $\Gamma_X$-acyclic resolution of the constant sheaf ${\bf R}$.

\hfill $\diamondsuit \diamondsuit \diamondsuit$

\medskip
Further, if we set $\Omega^0 = {\cal B}_{L, {\cal S}, X}$, and for brevity denote by

\bigskip
(4.3) \quad $ ( \Omega^*_{L, {\cal S}, X} ) $

\medskip
the complex in (4.1), then we can define the {\it cohomology algebra} $H^*_{L, {\cal S}, X}$ of (4.3). In this case,
the {\it abstract de Rham theorem}, see Mallios [1, chap iii, (3.25), or even (3.8)], becomes the relation

\bigskip
(4.4) \quad $ H^*(X,{\bf R}) ~=~ H^*_{de Rham}(X) ~=~ H^*_{L, {\cal S}, X} $

\medskip
which means that

\bigskip
(4.5) \quad $ \begin{array}{l}
                           the~usual~singular~or~\check{C}ech~cohomology~of~X,~as~well~as \\
                           the~standard~de~Rham~cohomology~of~X,~computed~in \\
                           terms~of~the~usual~smooth~functions~and~forms~on~X, \\
                           can~now~equally~be~computed~according~to~the~last,~and \\
                           highly~singular~term~in~(4.4).
                      \end{array} $

\medskip
That is, all the mentioned {\it cohomologies} of $X$, being actually functorially {\it isomorphic}, given that, by our
hypothesis, $X$ is paracompact and Hausdorff, see Mallios [1, chap. iii, (8.11)], they can now be computed through {\it
gene-\\ralized functions}, which can have {\it singularities} on arbitrary {\it dense} subsets in $X$, provided that the
complementary of the singularities is also dense in $X$. \\
\\
Now, the conclusion in (4.5) has implementations not only in mathematical physics, where it opens the way to dealing
with a large class of new and {\it enlarged} singularities. Indeed, it also has implications of a more purely mathematical
significance, due to its potential applicability in studying topology and differential geometry of several types of
nonsmooth spaces, and do so through cohomological methods of standard differential geometric character, see for
example Milnor classifying spaces, simplicial complexes, etc., and also Mostow, along with Rosinger [12], as well as
Mallios [1, chap. xi, sect. 12, in particular, (12.27)]. \\
\\
Similarly with Mallios \& Rosinger, with the consideration of the short exact exponential sheaf sequence associated
with (4.1) or (4.3), one can obtain further extensions of classical results in cohomology. \\
As an example of this type of results, one can refer to the following abelian {\it group isomorphism}, see also (2.38),
(2.41)

\bigskip
(4.6) \quad $ H^1 (X, ( {\cal C} ^\infty_X )^\bullet ) ~=~ H^2 (X, {\bf Z} ) ~=~ H^1 (X, {\cal B}^\bullet_{L, {\cal S}, X} ) $

\medskip
being of a purely {\it geometrical-physical} character. Indeed, in geometric terms (4.6) means that {\it any smooth}
${\bf R}$-{\it line bundle} on $X$ has a {\it coordinate 1-cocycle} consisting of elements locally from ${\cal B}_{L, {\cal
S}, X} )$. In physical terms, such as electromagnetism, for instance, the last abelian group in (4.6) classifies
cohomologically the {\it Maxwell fields}, within the abstract differential geometric setup of the present paper, which
involves the mentioned kind of new and enlarged singularities. More details on this, via a {\it gauge theoretic}
language, and as advocated in Mallios [1], can be found in Mallios [7,8].

\bigskip
{\bf Remark 2}

\medskip
It is important to note that, just like in Mallios \& Rosinger, where the nowhere dense differential algebras of
generalized functions were used, also in this paper, where the space-time foam differential algebras of generalized
functions are employed in the particular construction of the differential triad, there is again {\it no need} for any
topological algebra structure on the local sections $B_{L, {\cal S}} (U)$ of the sheaf ${\cal B}_{L, {\cal S}, X}$. This is
clearly unlike in the earlier formulations of the abstract differential geometric theory in Mallios [1-8], and it is, no doubt,
rather fortunate, since it further shifts the stress to purely algebraic ideas, concepts and constructs. \\
One of the reasons for the lack of need of any topological algebra structure on the algebras of generalized
functions under consideration is the following. It is becoming more and more clear that the classical Kuratowski-
Bourbaki topological concept is not suited to the mentioned algebras of generalized functions. Indeed, these
algebras prove to contain nonstandard type of elements, that is, elements which in a certain sense are infinitely
small, or on the contrary, large. And in such a case, just like in the much simpler case of nonstandard reals
$^*{\bf R}$, any topology which would be Hausdorff on the whole of the algebras of generalized functions, would by
necessity become discrete, therefore trivial, when restricted to usual, standard smooth functions, see for details
Biagioni. \\
\\
Here, in order to further clarify the issue of the possible limitations of the usual Kuratowski-Bourbaki concept of
topology, let us point out the following. Fundamental results from measure theory, predating the mentioned
concept of topology, yet having a clear topological nature, have never been given a suitable formulation within
that Kuratowski-Bourbaki concept. Indeed, such is the case, among others, with the Lebesgue dominated
convergence theorem, with the Lusin theorem on the approximation of measurable functions by continuous ones,
and with the Egorov theorem on the relation between pointwise and uniform convergence of sequences of
measurable functions. \\
Similar limitation of the Kuratowski-Bourbaki concept of topology appeared in the early 1950s, when attempts
were made to turn the convolution of Schwartz distributions into an operation simultaneously continuous in both
its arguments. More generally, it is well known that, given a locally convex topological vector space, if we
consider the natural bilinear form defined on its Cartesian product with its topological dual, then there will exist a
locally convex topology on this Cartesian product which will make the mentioned bilinear form simultaneously
continuous in both of its variables, if and only if our original locally convex topology is in fact as particular, as
being  a normed space topology. \\
It is also well known that in the theory of ordered spaces, in particular, ordered groups or vector spaces, there are
important concepts of convergence, completeness, boundedness, etc., which have never been given a suitable
formulation in terms of the Kuratowski-Bourbaki concept of topology. In fact, as seen in Oberguggenberger \&
Rosinger, powerful general results can be obtained about the existence of generalized solutions for very large
classes of nonlinear PDEs, by using order structures, and without any recourse to associated topologies. \\
Finally, it should be pointed out that, recently, differential calculus was given a new refoundation by using
standard concepts in category theory, such as naturalness. This approach also leads to topological type
processes, among them the so called toponomes or ${\cal C}$-spaces, which prove to be extensions of the usual
Kuratowski-Bourbaki concept of topology, see Nel, and the references cited there. \\
In this way, we can conclude that mathematics contains a variety of important {\it topological
type processes} which, so far, could not be formulated in convenient terms using the
Kuratowski-Bourbaki topological concept. And the differential algebras of generalized
functions, just as much as the far simpler nonstandard reals $^*{\bf R}$, happen to exhibit
such a class of topological type processes. \\
Unfortunately however, there seems so far to have been
insufficient awareness about the above state of affairs, a state which may be summarized as
follows :

\begin{itemize}

\item the Kuratowski-Bourbaki concept of topology is a rather narrow {\it particular} case of
the much large variety of {\it topological type processes} which have for long been used
successfully in various branches of mathematics,

\item important mathematical structures and developments cannot be confined within the limits
of the Kuratowski-Bourbaki concept of topology,

\item when using topological type structures beyond the Kuratowski-Bourbaki concept of
topology sufficient care should be taken in order to follow what may indeed be the {\it
naturally} extended concepts, thus avoiding to fall for one or another of the pet-concepts of
extended topology.

\end{itemize}

A recent systematic presentation of a wide range of topological type structures, together with
a number of their significant applications in Functional Analysis can be found in Beattie \&
Butzmann. One of the more important and useful extensions of Kuratowski-Bourbaki topology they
deal with is given by the so called {\it continuous convergence structures}, see their
Definition 1.1.5, on page 4. These topological type structures are specifically introduced in
order to deal with one of the long outstanding - even if less well known - major problems in
topology, namely, to define appropriate topological type structures on spaces of functions,
among them, on the spaces of continuous functions ${\cal C} ( X, Y )$, where $X$ and $Y$ are
topological spaces. \\
One of the main interests in dealing with that problem comes from the study of infinite
dimensional manifolds, where there has been a certain awareness related to the difficulties
coming from the limitation of the Kuratowski-Bourbaki concept of topology. In Kriegl \&
Michor, for instance, a significant effort was made in using a certain extended concept of
topology, called there a {\it convenient setting}, in order to deal with infinite dimensional
manifolds. As it happened however, this attempt is known not to have attained its ultimate
objectives. And one possible reason for that may precisely be in the insufficiently careful,
thus appropriate, or for that matter, convenient indeed, choice of the extended concept of
topology they happen to use. \\
On the other hand, the topological type processes on the nowhere dense differential algebras
of generalized functions, used in Mallios \& Rosinger, as well as on the space-time foam
differential algebras of generalized functions employed in this paper, can be given a suitable
formulation, and correspondingly, treatment, by noting that the mentioned algebras are in fact
{\it reduced powers}, see Bell \& Slomson, of ${\cal C}^\infty(X)$, and thus of ${\cal C}(X)$
as well. Let us give some further details related to this claim in the case of the space-time
foam algebras. The case of the nowhere dense algebras was treated in Mallios \& Rosinger. \\
Let us recall the definition in (2.12) of the space-time foam algebras, and note that it
obviously leads to

\bigskip
(4.7) \quad $ B_{L, {\cal S}} (X) ~=~ ({\cal C}^\infty(X))^\Lambda / {\cal J}_{L, {\cal S}} (X)  ~\subseteq~
                      ({\cal C}(X))^\Lambda / {\cal J}_{L, {\cal S}} (X) ~\subseteq~ \\
\\
 ~~~~~~~~~~~~~~~~~~~~~~~~~~~~~~~~~~~~~~~~~~~~~~~~~~~~~~~~~~~~~~~~~~ \subseteq~  {\cal C}(\Lambda \times X) / {\cal J}_{L, {\cal S}} (X) $

\medskip
assuming in the last term that on $\Lambda$ we consider the discrete topology. Now it is well known, Gillman \&
Jerison, that the algebra structure of ${\cal C}(\Lambda \times X)$ is connected to the topological structure of
$\Lambda \times X$, however, this connection is rather sophisticated, as essential aspects of it involve the
Stone-\u{C}ech compactification $\beta (\Lambda \times X)$ of $\Lambda \times X$. \\
It follows that a good deal of the discourse, and in particular, the topological type one, in the space-time foam algebras
$B_{L, {\cal S}} (X)$ may be captured by the topology of $\Lambda \times X$, and of course, by the far more involved
topology of $\beta(\Lambda \times X)$. Furthermore, the differential properties of these space-time foam algebras will,
in view of (2.18), (2.20), be reducible termwise to classical differentiation of sequences of smooth functions. \\
In short, in the case of the mentioned differential algebras of generalized functions, owing to their structure of
reduced powers, one obtains a 'two-way street'. Along it, on the one hand, the definitions and operations are
applied to sequences of smooth functions, and then reduced termwise to such functions, while on the other hand,
all that has to be done in a way which will be compatible with the 'reduction' of the 'power' by the quotient
constructions in (2.12), or in other words, (4.7). By the way, such a 'two-way street' approach has ever since the
1950s been fundamental in the branch of mathematical logic, called model theory, see Lo\u{s}. And a further quite
clear and detailed illustration of its workings can be seen in the next section, in the proof of Lemma 2. \\
But in order not to become unduly overwhelmed by ideas of model theory, let us recall here that the classical
Cauchy-Bolzano construction of the real numbers ${\bf R}$ is also a reduced power. Not to mention that a similar
kind of reduced power construction - in fact, its particular case called 'ultra-power' - gives the nonstandard reals
$^*{\bf R}$, as well.

\bigskip
{\bf Remark 3}

\medskip
The following lines of thought stand mainly in perspective with a potential application in the problem of quantization of
general relativity of the differential geometric framework that has been presented here, and whose {\it structure sheaf of
coefficients} was what we have called the sheaf of {\it space-time foam} algebras. \\
In this regard, by looking at the structure of the sheaf algebra ${\cal B}_{L, {\cal S}, X}$, see section 2.5, one obtains
by the definition of these algebras a decomposition of $X$ into singular and nonsingular domains, where $X$ is being
viewed as our space-time. And as seen, both of these domains can now be {\it dense} in $X$. \\
What is particularly important here is that such a decomposition is obtained by means of the same structure algebra
sheaf, which thus works simultaneously for both the regular and irregular parts of our space-time $X$. Furthermore, it is
worth noting that for both parts of $X$, our differential geometry, in other words, the corresponding {\it structure algebra
sheaf}, can equally be {\it commutative}, which may be referred to the earlier mentioned {\it Bohr's correspondence
principle}. What amounts to the same, and thus pertains to the same principle, is that one is in fact not compelled at all
to resort to a {\it noncommutative} structure algebra sheaf, as is traditionally done, in order to be able to cope with the
'quantum' part of $X$. Hence, one arrives at the mentioned principle which pertains to the {\it description of our
measurements} of a quantum system, and does so simply by the differential geometric apparatus presented here,
which at the end, it is but our algebra sheaf ${\cal B}_{L, {\cal S}, X}$. \\
In this connection we further remark that one can still formulate the corresponding {\it generalization of Einstein's
equations}, as this is fully explained in Mallios [7,8]. In fact, this can be done even within the so called {\it generalized
Lorentz} differential triad $(X, {\cal A}, \partial )$, see Mallios [6]. On the other hand, we still refer to Mallios [7,8],
concerning the initial formulation of the above, as well as for further pertinent comments. Other relevant argument can
be found in Heller [4] and Heller \& Sasin [4]. \\ \\

\bigskip
{\bf 5.~ Appendix}

\bigskip
{\bf Proof of Lemma 2}.

\medskip
It is easy to see that the restrictions in (2.37) will satisfy the required sheaf conditions. Indeed, what we actually prove
here is that the family (2.37) yields a {\it completepresheaf}, hence, equivalently, when categorically speaking, a sheaf,
see for instance, Mallios [1, chap. i, (11.37), or Theorem 13.1, and (13.18)]. \\
We can thus turn to check whether (2.37) satisfies the first of the two conditions on complete presheaves, related to
open covers, ( loc.cit., p. 46, Definition 11.1). \\
Let us therefore look at (2.41) and take $V = \bigcup_{i \in I}~V_i$, where $V_i \subseteq X$, with $i \in I$, are nonvoid
open. Given now two generalized functions $T,~ T^\prime \in B_{L, {\cal S}|_V}(V)$, let us assume that for $i \in I$, we
have

\bigskip
(5.1) \quad $ T|_{V_i} ~=~ T^\prime|_{V_i} $

\medskip
We then prove that

\bigskip
(5.2) \quad $ T ~=~ T^\prime $

\medskip
Let us note the relations

\bigskip
(5.3) \quad $ T = t + {\cal J}_{L, {\cal S}|_V}(V),~~~~~ T^\prime = t^\prime + {\cal J}_{L, {\cal S}|_V}(V) $

\medskip
with $t,~ t^\prime \in ({\cal C}^\infty(V))^\Lambda$, which follow from (2.12). Then (5.1) implies for $i \in I$

\bigskip
(5.4) \quad $ ( t^\prime - t ) |_{V_i} = w_i = ( w_{i ( \lambda )} ~|~ \lambda \in \Lambda ) \in {\cal J}_{L, {\cal S}|_{V_i}}(V_i) $

\medskip
Now for $i \in I$, we define the product mapping

\bigskip
(5.5) \quad $ {\cal C}^\infty_{V_i}(V) \times {\cal C}^\infty(V_i) \ni (\alpha,\psi) \longrightarrow
                                                                                         \alpha \psi \in {\cal C}^\infty(V) $

\medskip
where ${\cal C}^\infty_{V_i}(V)$ denotes the set of all smooth functions in ${\cal C}^\infty(V)$ whose support is in $V_i$,
while the product $\alpha \psi$ is defined by

$$ (\alpha \psi)(x) = \begin{array}{|l}
                                                          ~ \alpha(x) \psi(x) ~~~~~ \mbox{if}~ x \in V_i \\
                                                            \\
                                                          ~ 0 ~~~~~~~~~~~~~~~              \mbox{if}~ x \in V \setminus V_i
                                               \end{array} $$

\medskip
At this point, we consider a smooth partition of unity $( \alpha_l ~|~ l \in {\bf N} )$ on $V$, such that, see de Rham, each
$\alpha_l$ has a compact support contained in one of the $V_i$, and in addition, every point of $V$ has a neighbourhood which
intersects only a finite number of the supports of the various $\alpha_l$. In this way we obtain a mapping

\bigskip
(5.6) \quad $ N \ni l \longrightarrow i(l) \in I ~~~~\mbox{with}~~~~ \mbox{supp}~ \alpha_l \subseteq V_{i(l)} $

\medskip
Extending now termwise the above product to sequences of smooth functions indexed by $\lambda \in \Lambda$, we
can define the sequence of smooth functions, see (5.4)

\bigskip
(5.7) \quad $ w = ( w_\lambda ~|~ \lambda \in \Lambda ) = \sum_{l \in {\bf N}}~ \alpha_l~ w_{i(l)} \in
                                                                              ({\cal C}^\infty(V))^\Lambda $

\medskip
and then show that

\bigskip
(5.8) \quad $ w \in {\cal J}_{L, {\cal S}|_V}(V) $

\medskip
Once we have (5.8), we recall (5.4) - (5.6) and the fact that $( \alpha_l ~|~ l \in {\bf N} )$ is a smooth partition
of unity on $V$, and we obtain

$$ t^\prime - t = ( \sum_{l \in {\bf N}}~\alpha_l ) ( t^\prime - t ) = \sum_{l \in {\bf N}}~ \alpha_l
                 ( t^\prime - t ) = \sum_{l \in {\bf N}}~ \alpha_l ( t^\prime - t )|_{V_{i(l)}} = \sum_{l \in {\bf N}}~
                 \alpha_l w_{i(l)} = w $$

\medskip
which in view of (5.3) will indeed yield (5.2). \\
In order to obtain (5.8), and due to (2.10), (2.37), it suffices to find $\Sigma \in {\cal S}$, for which we have

\bigskip
(5.9) \quad $ w \in {\cal J}_{L, \Sigma \cap V}(V) $

\medskip
Let us therefore recall (5.4) and note that together with (5.5), (5.6) and (2.8), it results in

\bigskip
(5.10) \quad $ \alpha_l w_{i(l)} \in {\cal J}_{L, {\cal S}|_V}(V), ~~~\mbox{for}~~ l \in {\bf N} $

\medskip
thus (2.10) gives a sequence of sets $\Sigma^V_l \in {\cal S}|_V$, with $l \in {\bf N}$, such that

\bigskip
(5.11) \quad $ \alpha_l w_{i(l)} \in {\cal J}_{L, \Sigma^V_l}(V),~~~ \mbox{for}~~ l \in {\bf N} $

\medskip
and due to (5.10), (5.6), (2.8), we can further assume about $\Sigma^V_l$ that

\bigskip
(5.12) \quad $ \Sigma^V_l ~\subseteq~ \mbox{supp}~ \alpha_l, ~~~\mbox{for}~~ l\in {\bf N} $

\medskip
since we can always replace in (5.11) the initial $\Sigma^V_l$ with $\Sigma^V_l \cap \mbox{supp}~ \alpha_l$. However,
for $l \in {\bf N}$, we have $\Sigma^V_l = \Sigma_ l \cap V$, with suitable $\Sigma_l \in {\cal S}$. Then by taking

\bigskip
(5.13) \quad $ \Sigma = \bigcup_{l \in {\bf N}}~ \Sigma_l $

\medskip
and recalling (5.12), we obtain for $l \in {\bf N},~ x \in V,~ \Delta \subseteq V,~ \Delta$, neighbourhood of $x$, the
relations

$$ \Sigma_l \cap \Delta ~=~ ( \Sigma_l \cap V ) \cap \Delta ~=~ \Sigma^V_l \cap \Delta ~\subseteq~
                                                                                                                \mbox{supp}~ \alpha_l \cap \Delta $$

\medskip
It follows therefore from the assumed property of the supports of the partition of unity $( \alpha_l ~|~ l \in {\bf N} )$, that
for the given $V$, the sequence of singularity sets $\Sigma_l \in {\cal S}$, with $l \in {\bf N}$, satisfies condition (2.35).
Thus we have $\Sigma \cap V \in {\cal S}|_V$. But (5.13) yields

$$ \Sigma ~\cap~ V ~=~ \bigcup_{l \in {\bf N}}~ ( \Sigma_l ~\cap~ V ) ~=~ \bigcup_{l \in {\bf N}}~ \Sigma^V_l $$

\medskip
and thus (5.7),  (5.11) and (2.8) will give (5.9), and in this way, the proof of (5.8) is completed.

\bigskip
As a last step in order to show that (2.41) is a complete presheaf, let $T_i \in
B_{L, {\cal S}|_{V_i}}(V_i)$, with $i \in I$, be such that

\bigskip
(5.14) \quad $T_i|_{V_i \cap V_j} ~=~ T_j|_{V_i \cap V_j}$

\medskip
for all $i, j \in I$, for which $V_i \cap V_j \neq \phi$. Then we show that

\bigskip
(5.15) \quad $\begin{array}{l}
             \exists~ T \in B_{L, {\cal S}|_V}(V) ~: \\
             \\
             \forall~ i \in I ~: \\
             \\
             ~~~ T|_{V_i} ~=~ T_i
              \end{array}$

\medskip
Indeed, (2.12) gives the representations

\bigskip
(5.16) \quad $T_i ~=~ t_i + {\cal J}_{L, {\cal S}|_{V_i}}(V_i), ~~~\mbox{for}~~ i \in I$

\medskip
where $t_i \in ({\cal C}^\infty(V_i))^\Lambda$. But then (5.14) results in

\bigskip
(5.17) \quad $ ( t_i - t_j )|_{V_i \cap V_j} ~=~ w_{i~j} \in {\cal J}_{L, {\cal S}|_{V_i \cap V_j}}(V_i \cap V_j) $

\medskip
for all $i, j \in I$ such that $V_i \cap V_j \neq \phi$.

Let us take any fixed $i \in I$. Given $l \in {\bf N}$ such that $V_i \cap V_{i(l)} \neq \phi$, the relation (5.17) yields

$$ t_{i(l)} ~=~ t_i + w_{i(l)~i} ~~\mbox{on}~ V_i \cap V_{i(l)} $$

\medskip
thus (5.5), (5.6) lead to

$$ \alpha_l~t_{i(l)} ~=~ \alpha_l~t_i + \alpha_l~w_{i(l)~i} ~~\mbox{on}~ V_i $$

\medskip
But then

$$ \sum_{l \in {\bf N}}~ \alpha_l~t_{i(l)} ~=~ (\sum_{l \in {\bf N}}~ \alpha_l) t_i +
                                                   \sum_{l \in {\bf N}}~ \alpha_l~w_{i(l)~i} ~~\mbox{on}~ V_i $$

\medskip
or

\bigskip
(5.18) \quad $ \sum_{l \in {\bf N}}~ \alpha_l~t_{i(l)} ~=~ t_i + \sum_{l \in {\bf N}}~ \alpha_l~w_{i(l)~i}
                                                                                                     ~~\mbox{on}~ V_i $

\medskip
On the other hand, the relation

\bigskip
(5.19) \quad $ (\sum_{l \in {\bf N}}~ \alpha_l~w_{i(l)~i})|_{V_i} \in {\cal J}_{L. {\cal S}|_{V_i}}(V_i) $

\medskip
follows by an argument similar with the one we used for obtaining (5.8) via (5.9) - (5.13). In this way, if we define

\bigskip
(5.20) \quad $  T ~=~ t + {\cal J}_{L, {\cal S}|_V}(V) \in B_{L, {\cal S}_V}(V) $

\medskip
where

$$ t ~=~ \sum_{l \in {\bf N}}~ \alpha_l~ t_{i(l)} $$

\medskip
then (5.18) - (5.20) will give us (5.15), and the proof of the fact that (2.41) is a complete presheaf is completed.

\medskip
We turn now to proving that (2.41) is a {\it fine} sheaf. This however follows easily from (2.16), which as we have
noted, implies that $ 1_V = u(1) + {\cal J}_{L, {\cal S}|_V}(V)$ is the unit element in $B_{L, {\cal S}|_V}(V)$, thus
the partition of unity property together with (2.8) lead to

$$ 1_V ~=~ \sum_{l \in {\bf N}}~ (u(\alpha_l) + {\cal J}_{L, {\cal S}|_V}(V)) \in B_{L, {\cal S}|_V}(V) $$

\medskip
At this point, we are left only with showing that (2.41) is a {\it flabby} sheaf. Let $V^\prime \subseteq V \subseteq X$
be nonvoid open subsets, and let $T^\prime \in B_{L, {\cal S}|_{V^\prime}}(V^\prime)$. \\
Let us denote by $\Sigma^\prime$ the boundary of $V^\prime$ in $V$. Then clearly $\Sigma^\prime$ is closed and
nowhere dense in $V$, while $V \setminus ( V^\prime \cup \Sigma^\prime )$ is open in $V$. Further, since
$\Sigma^\prime$ is closed in $V$, there exists, Kahn, $\sigma^\prime \in {\cal C}^\infty(V)$, such that $\Sigma^\prime =
 \{~ x \in V ~|~ \sigma^\prime(x) = 0 ~\}$. \\
We shall use now an auxiliary function $\eta \in {\cal C}^\infty({\bf R})$ such that $\eta = 1$ on $(-\infty,-1] \cup [1,\infty)$,
while $\eta = 0$ on $[-1/2,1/2]$. And with its help, we can define the sequence of smooth functions $\beta_l \in {\cal C}^
\infty(V)$, with $l \in {\bf N}$, according to

\bigskip
(5.21) \quad $ \beta_l(x) ~=~  \begin{array}{|l}
                                                            ~~\eta( (l+1) \sigma^\prime(x) ) ~~~~\mbox{if}~~ x \in V^\prime \cup \Sigma^\prime \\
                                                            \\
                                                            ~~0 ~~~~~~~\mbox{if}~~ x \in V \setminus ( V^\prime \cup \Sigma^\prime )
                                              \end{array}$

\medskip
It is easy to check that

\bigskip
(5.22) \quad $ \mbox{supp}~ \beta_l ~\subseteq~ V^\prime,~~ \mbox{for}~ l \in {\bf N} $

\medskip
and

\bigskip
(5.23) \quad $ \begin{array}{l}
                                    \forall~ K \subset \subset V^\prime ~: \\
                                    \\
                                    \exists~ l \in {\bf N} ~: \\
                                    \\
                                    \forall~ k \in {\bf N},~ k \geq l ~: \\
                                    \\
                                    ~~~ \beta_k = 1 ~~~\mbox{on}~~ K
                         \end{array} $

\medskip
Let us now assume that $T^\prime$ has the representation, see (2.12)

\bigskip
(5.24) \quad $ T^\prime ~=~ t^\prime + {\cal J}_{L, {\cal S}|_{V^\prime}}(V^\prime) $

\medskip
where $t^\prime = (t^\prime_ \lambda ~|~ \lambda \in \Lambda) \in ({\cal C}^\infty(V^\prime))^\Lambda$. \\
Then we define

\bigskip
(5.25) \quad $ T ~=~ t + {\cal J}_{L, {\cal S}|_V}(V) \in B_{L, {\cal S}|_V}(V) $

\medskip
where $t = (t_\lambda ~|~ \lambda \in \Lambda)$ and, see (5.5), (5.22), $t_\lambda = \beta_{l_\lambda}~ t^\prime_
\lambda$, for $\lambda \in \Lambda$, $l_\lambda \in {\bf N}$ and $\lambda \leq l_\lambda$, this last inequality being
possible, since we assumed that ${\bf N}$ is cofinal in $\Lambda$.

\medskip
Before going further, we have to show that the definition of $T$ in (5.25) does not depend on the choice of $t^\prime$ in
(5.24). Let us therefore assume that, instead of the one in (2.35), we are given another representation

$$ T^\prime ~=~ t^{\prime *} + {\cal J}_{L, {\cal S}|_{V^\prime}}(V^\prime) $$

\medskip
with $t^{\prime *} = ( t^{\prime *}_\lambda ~|~ \lambda \in \Lambda ) \in ({\cal C}^\infty(V^\prime))^\Lambda$,
then clearly

\bigskip
(5.26) \quad $ t^{\prime *} - t^\prime = ( t^{\prime *}_\lambda - t^\prime_\lambda ~|~ \lambda \in \Lambda )
                                                                                                                    \in {\cal J}_{L, {\cal S}|_{V^\prime}}(V^\prime) $

\medskip
As above with $t$ in (5.25), let us now define $t^* = ( t^*_\lambda ~|~ \lambda \in \Lambda ) \in  ({\cal C}^\infty(V))^
\Lambda$ by $t^*_\lambda = \beta_{l_\lambda}~ t^{\prime *}_\lambda$, for $\lambda \in \Lambda$. \\
We show then that $T$ in (5.25) has also the representation

$$ T ~=~ t^* + {\cal J}_{L, {\cal S}|_V}(V) \in B_{L, {\cal S}|_V}(V) $$

\medskip
or that, equivalently

\bigskip
(5.27) \quad $ t^* -  t \in {\cal J}_{L, {\cal S}|_V}(V) $

\medskip
Indeed, from (5.26), (2.10), (2.37) we obtain $\Sigma \in {\cal S}$ such that

$$ t^{\prime *} - t^\prime \in {\cal J}_{L, \Sigma \cap V^\prime}(V^\prime) $$

\medskip
and then (2.8) gives

\bigskip
(5.28) \quad $ \begin{array}{l}
                               \forall~ x \in V^\prime \setminus ( \Sigma \cap V^\prime ) ~: \\
                               \\
                               \exists~ \lambda \in \Lambda ~: \\
                               \\
                               \forall~ \mu \in \Lambda,~ \mu \geq \lambda ~: \\
                               \\
                               \forall~ p \in {\bf N}^n ~: \\
                               \\
                               ~~~ D^p ( t^{\prime *}_\lambda - t^\prime_\lambda ) (x) = 0
                       \end{array} $

\medskip
But in view of (2.42), it follows that $( \Sigma \cap V ) \cup \Sigma^\prime \in {\cal S}|_V$, since clearly $\Sigma^\prime
\in {\cal S}_{nd}(V)$. Also, clearly, we have $t^*_\lambda - t_\lambda = \beta_{l_\lambda} ( t^{\prime *}_\lambda -
t^\prime_\lambda )$, with $\lambda \in \Lambda$. And then (5.28), (5.21) and (2.8) will directly lead to

$$ t^* - t \in {\cal J}_{L, ( \Sigma \cap V ) \cup \Sigma^\prime }(V) ~\subseteq~ {\cal J}_{L, {\cal S}|_V}(V) $$

\medskip
and the proof of (5.27) is completed.

\medskip
At last, it follows easily from (5.21) - (5.25) that

$$T|_{V^\prime} ~=~ T^\prime$$

\medskip
since a direct computation using also (2.8), gives

$$ t^\prime - t|_{V^\prime} \in {\cal J}_{L, \phi}(V^\prime) $$

\medskip
In this way the flabbiness of (2.41) is proved.

\end{document}